\documentclass[leqno,draft]{article}

\def\dist{\mathop{\rm dist}}

\begin{document}

\title{A geometric construction for fractal sets and related linear
operators}

\author{Stephen Semmes \\
        Rice University}

\date{}

\maketitle

        Let $U$ be a nonempty bounded open set in ${\bf R}^n$, and put
$E = \partial U$.  Thus $E$ is a nonempty compact set in ${\bf R}^n$
with empty interior.  By standard arguments,
\begin{equation}
        d_E(x) = \dist(x, E) = \inf \{|x - y| : y \in E\}
\end{equation}
is a Lipschitz function on ${\bf R}^n$ of order $1$ and with constant $1$.
Of course,
\begin{equation}
        \dist(x, {\bf R}^n \backslash U) = \dist(x, E)
\end{equation}
when $x \in U$.  If $V$ is a connected component of $U$, then
\begin{equation}
        \dist(x, {\bf R}^n \backslash V) = \dist(x, \partial V) = \dist(x, E)
\end{equation}
when $x \in V$.

        Let us identify ${\bf R}^{n + 1}$ with ${\bf R}^n \times {\bf
R}$, so that an element of ${\bf R}^{n + 1}$ may be expressed as $(x,
t)$ with $x \in {\bf R}^n$ and $t \in {\bf R}$.  Using these
coodinates, we can also identify ${\bf R}^n$ with the $t = 0$
hyperplane in ${\bf R}^{n + 1}$.  Put
\begin{equation}
        \widehat{U} = \{(x, t) \in {\bf R}^{n + 1} : x \in U, \ |t| < d_E(x)\}.
\end{equation}
Similarly, put
\begin{equation}
        \widehat{V} = \{(x, t) \in {\bf R}^{n + 1} : x \in V, \ |t| < d_E(x)\}
\end{equation}
for every connected component $V$ of $U$.

        Clearly $\widehat{U}$ is a bounded open set in ${\bf R}^{n +
1}$.  For each connected component $V$ of $U$, $\widehat{V} \subseteq
\widehat{U}$ is a connected open set in ${\bf R}^{n + 1}$, since $V$
is.  If $V$, $W$ are distinct connected components of $U$, then $V
\cap W = \emptyset$, and hence $\widehat{V} \cap \widehat{W} =
\emptyset$.  Using the fact that
\begin{equation}
        U = \bigcup \{V : V \hbox{ is a connected component of } U\},
\end{equation}
one can check that
\begin{equation}
        \widehat{U} = \bigcup \{\widehat{V} :
                                V \hbox{ is a connected component of } U \}.
\end{equation}
It follows that the connected components of $\widehat{U}$ are the sets
$\widehat{V}$ corresponding to the connected components $V$ of $U$.

        Let $\widetilde{E}$ be the boundary of $\widehat{U}$ in ${\bf
R}^{n + 1}$, a nonempty compact set with empty interior in ${\bf R}^{n
+ 1}$.  Thus $\widetilde{E}$ consists of the $(x, t) \in {\bf R}^{n +
1}$ such that $x$ is in the closure of $U$ and $t = \pm d_E(x)$, and
the intersection of $\widetilde{E}$ with the $t = 0$ hyperplane is
equal to $E$.  As a simple variant of this construction, one might
take $d_E(x)$ to be a regularized version of the distance to $E$, so
that the boundary of $\widehat{U}$ is smooth away from the $t = 0$
hyperplane.  One might also prefer to choose the components
$\widehat{V}$ of $\widehat{U}$ to have some specific form, depending
on the situation.  If a component $V$ of $U$ is an $n$-dimensional
topological ball, then the corresponding component $\widehat{V}$ of
$\widehat{U}$ is an $(n + 1)$-dimensional topological ball.  In any
case, the boundary of $\widehat{V}$ is the same topologically as
doubling the closure of $V$, which is to say gluing two copies of the
closure of $V$ together along the boundary.  If $E$ is the Cantor set
and $U$ is its complement in the unit interval, then $\widetilde{E}$
is connected but still somewhat complicated.  Note that $E$ has
topological dimension $n - 1$ under the conditions mentioned at the
beginning, and $\widetilde{E}$ has dimension $n$.

        One can consider Toeplitz-type operators associated to
$\widehat{U}$ and $\widetilde{E}$, using quaternionic or Clifford
analysis.  One could have done this already with $U$ and $E$, but it
is sometimes of interest to change dimensions, and switch between even
and odd dimensions in particular.  There are also subtleties in
general about exactly what functions are used to define Toeplitz and
related operators, and which functions are in the domains of these
operators and which are in the ranges.

\end{document}